\documentclass[a4paper,11pt]{article}
\usepackage{amsfonts,amssymb,amsmath, dsfont}
\usepackage{bm}
\usepackage{color}
\usepackage{upgreek}
\usepackage[english]{babel}
\usepackage{graphicx}
 \makeatletter
 \@addtoreset{equation}{section}
 \makeatother

 \makeatletter
 \@addtoreset{enunciato}{section}
 \makeatother

 \newcounter{enunciato}[section]

 \newtheorem{ittheorem}{Theorem}
 \newtheorem{itlemma}{Lemma}
 \newtheorem{itproposition}{Proposition}
 \newtheorem{itcorollary}{Corollary}
 \newtheorem{itdefinition}{Definition}
 \newtheorem{itremark}{Remark}
 \newtheorem{itclaim}{Claim}
 \newtheorem{itfact}{Fact}
 \newtheorem{itconjecture}{Conjecture}

 \newenvironment{theorem}{\addtocounter{enunciato}{1}
 \begin{ittheorem}}{\end{ittheorem}}

 \newenvironment{lemma}{\addtocounter{enunciato}{1}
 \begin{itlemma}}{\end{itlemma}}

 \newenvironment{proposition}{\addtocounter{enunciato}{1}
 \begin{itproposition}}{\end{itproposition}}

 \newenvironment{corollary}{\addtocounter{enunciato}{1}
 \begin{itcorollary}}{\end{itcorollary}}

 \newenvironment{definition}{\addtocounter{enunciato}{1}
 \begin{itdefinition}}{\end{itdefinition}}

 \newenvironment{remark}{\addtocounter{enunciato}{1}
 \begin{itremark}}{\end{itremark}}

 \newenvironment{claim}{\addtocounter{enunciato}{1}
 \begin{itclaim}}{\end{itclaim}}

 \newenvironment{fact}{\addtocounter{enunciato}{1}
 \begin{itfact}}{\end{itfact}}

 \newenvironment{conjecture}{\addtocounter{enunciato}{1}
 \begin{itconjecture}}{\end{itconjecture}}

 \newcommand{\be}[1]{\begin{equation}\label{#1}}
 \newcommand{\ee}{\end{equation}}

 \newcommand{\bl}[1]{\begin{lemma}\label{#1}}
 \newcommand{\el}{\end{lemma}}

 \newcommand{\br}[1]{\begin{remark}\label{#1}}
 \newcommand{\er}{\end{remark}}

 \newcommand{\bt}[1]{\begin{theorem}\label{#1}}
 \newcommand{\et}{\end{theorem}}

 \newcommand{\bd}[1]{\begin{definition}\label{#1}}
 \newcommand{\ed}{\end{definition}}

 \newcommand{\bcl}[1]{\begin{claim}\label{#1}}
 \newcommand{\ecl}{\end{claim}}

 \newcommand{\bfact}[1]{\begin{fact}\label{#1}}
 \newcommand{\efact}{\end{fact}}

 \newcommand{\bp}[1]{\begin{proposition}\label{#1}}
 \newcommand{\ep}{\end{proposition}}

 \newcommand{\bc}[1]{\begin{corollary}\label{#1}}
 \newcommand{\ec}{\end{corollary}}

 \newcommand{\bcj}[1]{\begin{conjecture}\label{#1}}
 \newcommand{\ecj}{\end{conjecture}}

 \newcommand{\bpr}{\begin{proof}}
 \newcommand{\epr}{\end{proof}}

 \newcommand{\bprlem}[1]{\begin{proofof}{\it Lemma \ref{#1}}.\,\,}
 \newcommand{\eprlem}{\end{proofof}}

 \newcommand{\bprthm}[1]{\begin{proofof}{\it Theorem \ref{#1}}.\,\,}
 \newcommand{\eprthm}{\end{proofof}}

 \newcommand{\bprprop}[1]{\begin{proofof}{\it Proposition \ref{#1}}.\,\,}
 \newcommand{\eprprop}{\end{proofof}}

 \newcommand{\bi}{\begin{itemize}}
 \newcommand{\ei}{\end{itemize}}

 \newcommand{\ben}{\begin{enumerate}}
 \newcommand{\een}{\end{enumerate}}

 \newenvironment{proof}{\noindent {\em Proof}.\,\,}{\hspace*{\fill}$\halmos$\medskip}
 \newenvironment{proofof}{\noindent {\em Proof of\,\,}}{\hspace*{\fill}$\halmos$\medskip}
 \newcommand{\halmos}{\rule{1ex}{1.4ex}}

 \newcommand{\one}{{\mathchoice {1\mskip-4mu\mathrm l}
         {1\mskip-4mu\mathrm l}
         {1\mskip-4.5mu\mathrm l}
         {1\mskip-5mu\mathrm l}}}

\def \E {{\mathbb E}}

\def \N {{\mathbb N}}
\def \P {{\mathbb P}}

\def \ba {\begin{array}}
\def \ea {\end{array}}

\def \T {{\mathbb{T}}}

\def\one{\rlap{\mbox{\small\rm 1}}\kern.15em 1}

\begin{document}
\title{Supercriticality of Annealed Approximations of Boolean Networks}

\author{Thomas Mountford\textsuperscript{1} and Daniel Valesin\textsuperscript{2}}

\footnotetext[1]{\'Ecole Polytechnique F\'ed\'erale de Lausanne,
D\'epartement de Math\'ematiques,
1015 Lausanne, Switzerland}
\footnotetext[2]{University of British Columbia, Department of Mathematics, V6T1Z2 Vancouver, Canada}
\date{November 22, 2012}
\maketitle

\begin{abstract}
We consider a model recently proposed by Chatterjee and Durrett \cite{CD} as an ``annealed approximation'' of boolean networks, which are a class of cellular automata on a random graph, as defined by S. Kauffman \cite{kauff}. The starting point is a random directed graph on $n$ vertices; each vertex has $r$ input vertices pointing to it. For the model of \cite{CD}, a discrete time threshold contact process is then considered on this graph: at each instant, each vertex has probability $q$ of choosing to receive input; if it does, and if at least one of its input vertices were in state 1 at the previous instant, then it is labelled with a 1; in all other cases, it is labelled with a 0. $r$ and $q$ are kept fixed and $n$ is taken to infinity. Improving a result of \cite{CD}, we show that if $qr > 1$, then the time of persistence of activity of the dynamics is exponential in $n$. 
\end{abstract}

{\bf\large{}}\bigskip


\section{Introduction}
Random boolean networks were introduced by Stuart Kauffman in 1969 \cite{kauff} as models of gene regulatory networks. A gene regulatory network is a set of genes in a cell that iteratively communicate with each other, using their RNA transcripts as messages, and this communication affects each gene's activity. They are thus information networks and control systems for the activity of the cell.

Let us define Kauffman's model. The following definition depends on three parameters: $n,\; r \in \mathbb{N}$ with $r \leq n$ and $p \in (0,1)$ (though Kauffman only considered the case $p = 1/2$). The letters $a, b$ will denote two possible states of a gene. Let $V_n = \{x_1, \ldots, x_n\}$ be the set of genes. For each $x \in V_n$, we independently choose:\medskip\\
$\bullet\;$ a set $y(x) = \{y_1(x), \ldots, y_r(x)\} \subset V_n -\{x\}.$ The choice is made uniformly among all possibilities. $y(x)$ is called the \textit{influence set} of $x$.\medskip\\
$\bullet\;$ a function $f_x: \{a,b\}^r \to \{a,b\}.$ The values $\{f_x(\omega): \omega \in \{a,b\}^{y(x)}\}$ are chosen independently, with probability $p$ to be equal to $a$ and $1-p$ to be equal to $b$.\medskip\\
Having made all these random choices, and given an initial configuration $\eta_0 \in \{a,b\}^{V_n}$, we define a deterministic, discrete time dynamics $(\eta_t)_{t = 0, 1,\ldots}$ in $\{a,b\}^{V_n}$ by putting
$$\eta_{t+1}(x) = f_x\big(\eta_t(y_1(x)), \ldots, \eta_t(y_r(x))\big), \qquad t \geq 0.$$
That is, at each instant we verify the previous states in the influence set of $x$ and from these, determine the state of $x$ using the function $f_x$.  

Since the evolution is deterministic and the state space is finite, every initial configuration is in the domain of attraction of a periodic orbit or a fixed point. Typical aspects of interest in random boolean networks are the number of these attractors, their stability, periods and the time to reach them. As thoroughly explained in \cite{kauffman93}, simulations of the model suggested the existence of two regimes, depending on the choice of parameters, in which drastically different behaviours arise: in the \textit{ordered} (or subcritical) regime, the orbits and the typical time to reach them grow slowly with $n$, whereas in the \textit{chaotic} (or supercritical) regime, they grow rapidly with $n$.

In \cite{derpo}, Derrida and Pomeau proposed an ``annealed approximation'' of random boolean networks; in it, the random aspects of the network (namely, the underlying graph and the rules of evolution) are updated at each time step instead of remaining fixed. This simplification destroys important correlations in the system, but still allows for a rigorous proof of a phase transition given by a curve that agrees with simulations, $2rp(1-p) = 1$ (the ordered regime corresponding to $2rp(1-p) < 1$). Chatterjee and Durrett proposed in \cite{CD} a new ``annealed approximation'' model and proved the phase transition with the same critical curve. Their model was a more accurate approximation because, though the rules of evolution were resampled with time, the random graph was kept fixed. The resulting dynamics was that of a threshold contact process on the random graph, and allowed for new insight by providing an analogy between the flow of information in random boolean networks and the evolution of branching processes.

We now define the model of \cite{CD}. In a short subsection at the end of the Introduction, we clarify the precise connection between this model and Kauffman's original boolean networks. We start with parameters $n,\;r \in \mathbb{N}$ with $r \leq n$ and $q \in (0,1)$ (in the comparison with boolean networks, $q$ plays the role of $2p(1-p)$). Define the random graph $G_n = (V_n, E_n)$ exactly as before. We will now define a discrete time Markov chain $(\xi_t)_{t \geq 0}$ with state space $\{0,1\}^{V_n}$ and initial configuration $\xi_0 \equiv 1$. Its transition kernel is given by
$$P(\xi, \xi') = \left(\prod_{\substack{x \in V_n:\sum_i \xi(y_i(x)) = 0}}I_{\{\xi'(x) = 0\}} \right)\left(\prod_{\substack{x \in V_n:\sum_i \xi(y_i(x)) > 0}}\left( q \cdot I_{\{\xi'(x) = 1\}} + (1-q) \cdot I_{\{\xi'(x) = 0\}}\right) \right),$$
where $\xi, \xi' \in \{0,1\}^{V_n}$ and $I$ denotes the indicator function. It will be useful to construct this Markov chain with a set of auxiliary Bernoulli random variables. Let $\{B^x_t: x \in V_n,\;t\geq 1\}$ be a family of independent Bernoulli random variables with parameter $q$; given $\xi_t \in \{0,1\}^{V_n}$, we put
$$\xi_{t+1}(x) = \left\{\begin{array}{ll}1&\text{if } B^x_{t+1} = 1 \text{ and } \sum_{i=1}^r \xi_t(y_i(x))>0;\\0&\text{otherwise}. \end{array} \right.$$
When $B^x_t = 1$, we say that $x$ \textit{receives input} at time $t$; therefore, a vertex is set to 1 if and only if it receives input at that time and at least one of its input vertices $y_1(x),\ldots, y_r(x)$ was set to 1 at the previous time. $\P_n$ will denote a probability measure both for the choice of $G_n$ and for the family $\{B^x_t\}$ (they are of course taken independently). We sometimes abuse notation and associate $\xi \in \{0,1\}^{V_n}$ with $\{x\in V_n: \xi(x) = 1\}.$

It is readily seen that the identically zero configuration is absorbing for this chain and that it is eventually reached with probability 1. In \cite{CD}, the authors study the time $\uptau_n$ it takes for this to occur and the typical proportion of sites that are in state 1 at times before $\uptau_n$. By a simple comparison between the time dual of the model (as defined below) and a subcritical branching process, it is easy to show that, if $qr < 1$, then $\uptau_n$ behaves as $\log n$, and this is associated to the ordered regime of random boolean networks. In \cite{CD}, the following result is shown, characterizing the chaotic regime. Let $\rho = \rho(q, r)$ denote the probability of survival for a branching process in which individuals have probability $q$ of having $r$ children and probability $1-q$ of having none. Let $|A|$ denote the cardinality of the set $A$.\medskip

\noindent \textbf{Theorem.} \cite{CD} \textit{If $q(r-1) > 1$, then for every $\epsilon > 0$ there exists $c > 0$ such that, as $n \to \infty$},
$$\inf_{0 \leq t \leq e^{cn}} \P_n \left(\frac{|\xi_t|}{n} \geq \rho - \epsilon \right) \stackrel{n \to \infty}{\xrightarrow{\hspace*{0.8cm}}} 1.$$
Under the more general hypothesis $qr > 1$, only a weaker result was obtained in \cite{CD} because of certain technical difficulties related to the structure of the graph $G_n$ and the comparison to the branching process. We have dealt with these and obtained the stronger result:
\begin{theorem}
\label{thm:main}
If $qr > 1$, then there exists $c > 0$ such that, for any $\epsilon > 0$ and any sequence $(t_n)$ with $t_n \to \infty$ and $t_n \leq e^{cn}$,
$$\inf_{t_n\leq t \leq e^{cn}}\P_n\left(\rho - \epsilon < \frac{|\xi_t|}{n} < \rho+\epsilon \right) \stackrel{n \to \infty}{\xrightarrow{\hspace*{0.8cm}}} 1.$$
\end{theorem}

To explain why this result is to be expected and, in particular, the link with the mentioned branching process, we introduce the time dual of the process. Fix a realization of $G_n = (V_n, E_n)$ and $\{B^x_t: x\in V_n,\;t\geq 1\}$, define $\hat E_n$ as the set of directed edges obtained by inverting the edges of $E_n$ and $\hat G_n = (V_n, \hat E_n)$. Note that 
$$\{y_i(x): 1 \leq i \leq r\} = \{z: (x, z) \in \hat E_n\};$$
that is, in $\hat G_n$ each vertex ``points to'' $r$ vertices. Fix $T > 0$ and put $\hat B^{x,T}_t = B^x_{T-t}$ for $0 \leq t < T$. Given $A \subset V_n$, define $\hat \xi^{A,T}_0 = I_A$ and, for $0 \leq t < T$,
$$\hat \xi^{A,T}_{t+1}(z) = \left\{\begin{array}{ll}1&\text{if for some $x, i$, we have $y_i(x) = z,\;\hat\xi^{A,T}_t(x)=1$ and $\hat B^{x,T}_t = 1;$}\\0&\text{otherwise}. \end{array} \right.$$
When $\hat \xi^{A,T}_t(x) = 1$ and $\hat B^{x,T}_t = 1$, we say that $x$ \textit{gives birth} at time $t$. Let us describe the dual dynamics in words. Given the configuration $\hat \xi_t$, we go over every vertex that is in state 1 and determine which of them give birth at time $t$ -- for each vertex, this happens with probability $q$ and independently. For each vertex $x$ that gives birth at time $t$, we set the vertices $y_1(x),\ldots, y_r(x)$ to 1 at time $t+1$. Vertices that are not set to 1 by this procedure are then set to 0.

We have the \textit{duality equation}
$$\left\{\xi_T \cap A \neq \varnothing\right\} = \left\{\hat \xi^{A,T}_T \neq \varnothing\right\},$$
a consequence of which is that, under $\P_n$, $|\xi_T|$ and $|\{x: \hat \xi^{\{x\},T}_T \neq \varnothing\}|$ have the same distribution. 
 
Since we will mostly work with the dual process and will rarely have to consider the primal and dual processes jointly, we drop the superscript $T$ and assume that $\hat \xi^A_t$ is defined for all positive times with the evolution rule explained above. If $A = \{x\}$, we write $\hat \xi^x_t$. $\P_n$ will still denote a probability measure for both the random graph $\hat G_n$ and the dual process. For a fixed realization $\hat G_n$ of the graph, we will also need the \textit{quenched measure} $P_{\hat G_n}$, under which the dual process with specified initial configuration will be defined on this graph.

Now, assume that $n$ is very large with respect to $r$. If $g$ is another integer that is much larger than $r$ and much smaller than $n$, then with high probability, the subgraph of $\hat G_n$ with vertex set
$$\{z \in V_n: \text{for some }k\leq g \text{ and } z_1, \ldots, z_k \in V_n, \text{ we have } x \to z_1 \to \cdots \to z_k \to z \text{ in } \hat G_n\}$$
and edge set equal to the set of edges of $\hat E_n$ that start and end at vertices in the above set will simply be a directed tree of degree $r$ rooted in $x$. Conditioning on the event that this subgraph is indeed a tree, the evolution of $|\hat \xi^x_t|$ up to time $g$ will be exactly that of the branching process mentioned before Theorem \ref{thm:main}. In addition, it is not difficult to see that, without any conditioning, $|\hat \xi^x_t|$ is stochastically dominated by such a process. These remarks clarify why the model exhibits two phases in exact correspondence with the branching process. If the expected offspring size $qr < 1$, then $\hat \xi^x_t$ dies out faster than the corresponding subcritical branching process, and the primal $\xi_t$ rapidly reaches the zero state. On the other hand, if $qr>1$, the above theorem states that the system survives for a time that is exponentially large in $n$, characterizing the supercritical regime.

The structure of our proof is similar to that of \cite{CD}. First, using the comparison with the branching process and a second moment argument, we show that with probability tending to 1 as $n \to \infty$, the set of vertices $S = \{x:|\hat \xi^x_t|\text{ eventually reaches }(\log n)^2\}$ has size close to $\rho\cdot n$ (Proposition \ref{prop:bfsn}). Second, we show that with probability tending to 1 as $n \to \infty$, the graph $\hat G_n$ is such that, for any $A \subset V_n$ with $|A| \geq (\log n)^2$, the probability that $\hat \xi^A_t$ remains active up to time $e^{cn}$, for some fixed constant $c$, is larger than $n^{-\sqrt{\log n}}$ (Proposition \ref{prop:afsn}). We can then use a simple union bound to argue that with high probability, for every $x$ in $S$, $(\hat \xi^x_t)$ remains active until $e^{cn}$, and conclude by duality.

Our main contribution is Proposition \ref{prop:afsn}; let us briefly explain the ideas that go into its proof. Given $A \subset V_n$, suppose we reveal, one by one, the elements of the set $A_1 = \{y_i(x): 1\leq i \leq r,\;x\in A\}$, then $A_2 = \{y_i(x): 1\leq i \leq r,\;x\in A_1\}$, until $A_g$, for some fixed $g \in \N$. Let $B(A, g)$ be the subgraph of $\hat G_n$ with vertex set $A \cup A_1 \cup \cdots \cup A_g$ and edge set equal to the edges of $\hat E_n$ which start and end at vertices in this set. For most choices of $A$, $B(A,g)$ is just a disjoint union of $|A|$ directed trees, so that $\{|\hat \xi^A_t|\}_{0 \leq t \leq g}$ is exactly a branching process. However, for some choices of $A$, when revealing $A_1, \cdots, A_g$, we will see some ``collisions'', that is, some vertices will be found more than once. We say that $A$ is \textit{expansive} if the number of collisions is not too large, so that $\{|\hat \xi^A_t|\}_{0 \leq t \leq g}$ is not too far from the branching process and consequently, $|\hat \xi^A_g|$ is very likely to be larger than $|A|$ (see Lemma \ref{lem:robgood}). We then show that, with high probability, for some $c > 0$, there is no set $A \subset V_n$ with $(\log n)^2 \leq |A| \leq cn$ that is not expansive (Lemma \ref{lem:robeasy}). It is then quite easy to put Lemmas \ref{lem:robgood} and \ref{lem:robeasy} together to obtain Proposition \ref{prop:afsn}.

\subsection{Relationship between the boolean network model and the threshold contact process}
Assume that $(\eta_t)_{t \geq 0}$ is defined as in the beginning of the Introduction and the initial configuration is chosen with the product measure $\otimes_{x \in V_n} (p \cdot \delta_{\{a\}} + (1-p) \cdot \delta_{\{b\}})$. We will now show how, through a sequence of simplifications, this model is reduced to the threshold contact process $(\xi_t)_{t \geq 0}$.
\begin{enumerate}
\item[1)] Instead of maintaining the functions $\{f_x: x \in V_n\}$ fixed from the start, we can define a random dynamics $(\eta^2_t)$ with the property that each value $f_x(\omega)$ is resampled each time it is queried. That is, we select the influence sets and an initial condition $\eta_0$ as before, and then define
$$\begin{aligned}
\qquad\qquad&\eta^2_0 \equiv \eta_0;\\
&\eta^2_{t+1}(x) = \left\{\begin{array}{ll}\eta^2_t(x) &\text{if } \eta^2_t(y_i(x)) = \eta^2_{t-1}(y_i(x)) \; \forall i;\medskip\\ \begin{array}{c}a \text{ with probability } p,\\ b \text{ with probability } 1-p \end{array}&\text{ otherwise.}\end{array}\right.
\end{aligned}$$
Of course, this simplification will have a smaller effect the larger $r$ is, because then we will have fewer repetitions of $(\eta_t(y_1(x)),\ldots, \eta_t(y_r(x)))$ as $t$ varies.

\item[2)] With the above definition, the dynamics of $(\eta^2_t)$ becomes constant after an instant $t$ such that $\eta^2_t(x) = \eta^2_{t-1}(x)$ for all $x$, and this almost surely happens in finite time. In any case, we will have a decrease or interruption in the \textit{activity} of the system, meaning that the state of most sites (or all of them) will stay unchanged from one period to the other. If we only want to study the time it takes for activity to decrease or to halt, it suffices that we consider an auxiliary process $(\Gamma \eta^2_t)_{t \geq 0}$ in $\{0,1\}^{V_n}$ which marks whether or not the state of each vertex has changed from one period to the other. Formally, given $(\eta^2_t)_{t \geq 0}$, put $\Gamma \eta^2_0 \equiv 1$ and $\Gamma \eta^2_{t+1}(x) = I_{\{\eta^2_{t+1}(x) \neq \eta^2_t(x)\}}$, where $I$ denotes the indicator function. 

\item[3)] Condition on the choice of the influence sets $\{y(x): x \in V_n\}$. The process $\left(\eta^2_t,\; \Gamma \eta^2_t \right)_{t \geq 0}$ is then a Markov chain, but neither $\left(\eta^2_t \right)_{t \geq 0}$ nor $\left(\Gamma \eta^2_t \right)_{t \geq 0}$ in isolation is Markovian. Our last step will be simplifying $(\Gamma \eta^2_t)$ so that it becomes a Markov chain $(\Gamma^2_t)$. To this end, note that, for any $z$,
$$\begin{aligned}
&\P\left(\Gamma \eta^2_{t+1}(z) = 1 \; {\big \vert}\; (y(x))_{x \in V_n}, \; \Gamma \eta^2_t \right) \\[+5pt]
&= \left\{
\begin{array}{ll}
0& \text{if } \Gamma \eta^2_t(y_i(z)) = 0 \; \forall i;\\[+5pt]
\begin{array}{l}\P(\;\eta^2_{t+1}(z) = a,\; \eta^2_t(z) = b\; {\big \vert}\; (y(x))_{x \in V_n}, \; \Gamma \eta^2_t\;) \\+\;\P(\;\eta^2_{t+1}(z) = b,\; \eta^2_t(z) = a\; {\big \vert}\; (y(x))_{x \in V_n}, \; \Gamma \eta^2_t\;)
\end{array}&\text{otherwise}.
\end{array} \right.\\[+5pt]
&= \left\{
\begin{array}{ll}
0& \text{if } \Gamma \eta^2_t(y_i(z)) = 0 \; \forall i;\\[+5pt]
\begin{array}{l}p\cdot \P(\; \eta^2_t(z) = b\; {\big \vert}\; (y(x))_{x \in V_n}, \; \Gamma \eta^2_t\;) \\+\;(1-p)\cdot\P(\; \eta^2_t(z) = a\; {\big \vert}\; (y(x))_{x \in V_n}, \; \Gamma \eta^2_t\;)
\end{array}&\text{otherwise}.
\end{array} \right.
\end{aligned}$$
The idea of defining $(\Gamma^2_t)$ is using the above expression but pretending that $\P(\; \eta^2_t(z) = a\; {\big \vert}\; (y(x))_{x \in V_n}, \; \Gamma \eta^2_t\;)$ and $\P(\;\eta^2_t(z) = b\; {\big \vert}\; (y(x))_{x \in V_n}, \; \Gamma \eta^2_t\;)$ are always equal to their values when $t=0$: $p$ and $1-p$, respectively. We thus get $\Gamma^2_0 \equiv 1$ and
$$\begin{aligned}
&\P\left(\;\Gamma^2_{t+1}(z) = 1 \; {\big \vert} \; (y(x))_{x \in V_n},\;\Gamma^2_t\;\right) = \left\{\begin{array}{ll}0&\text{if } \Gamma^2_t((y_i(z))=0\; \forall i;\\[+5pt] 2p(1-p)&\text{otherwise.}  \end{array}\right.
\end{aligned}$$
\end{enumerate}
Then, $(\Gamma^2_t)_{t \geq 0}$ with parameters $r, p$ has the law of the process $(\xi_t)_{t \geq 0}$ with parameters $r, q = 2p(1-p)$. This process is thus a simplified, randomly-evolving model for the activity of boolean networks.

When $p = 1/2$, point $3)$ has no effect: $(\Gamma \eta^2_t)$ and $(\Gamma^2_t)$ have the same distribution because, by symmetry between states $a$ and $b$, the condition $$\P(\; \eta^2_t(z) = a\; {\big \vert}\; (y(x))_{x \in V_n}, \; \Gamma \eta^2_t\;) = \P(\; \eta^2_t(z) = a\; {\big \vert}\; (y(x))_{x \in V_n}, \; \Gamma \eta^2_t\;) = 1/2$$
holds true in this case.
When $p$ is taken far from $1/2$, point $3)$ may lead to a worse approximation, because there is no reason to assume that the conditional distribution of $\eta^2_t(z)$ remains near its starting point.

\section{Proof of Theorem \ref{thm:main}}
In all results and proofs that follow, we assume that $qr>1$. We will write
$$k_n = (\log n)^2,\qquad s_n = (\log \log n)^2.$$
We start with two propositions that together will yield Theorem \ref{thm:main}. Proposition \ref{prop:bfsn} is proved essentially by a repetition of arguments in \cite{CD}; we include a brief proof for completeness.

\begin{proposition} For any $\epsilon > 0$,
$$\lim_{n \to \infty} \P_n\left(\rho-\epsilon < \frac{\left|\left\{x \in V_n: |\hat \xi^x_{s_n}| > k_n \right\}\right|}{n}  < \rho + \epsilon \right) = 1.$$
\label{prop:bfsn}
\end{proposition}

\begin{proof}
Let $(Z_t)_{t = 0, 1,\ldots}$ be the branching process with $Z_0 = 1$ and the offspring distribution that gives mass $q$ to $r$ and $1-q$ to $0$, so that $\rho = \P(Z_t \neq 0 \; \forall t) > 0$. On the event $\{Z_t \neq 0 \;\forall t\}$, $\frac{Z_t}{(qr)^t}$ almost surely converges to a positive limit (see \cite{durprob}, Section 5.3.4). Thus, defining $\rho_n = \P(|Z_{s_n}| > k_n)$ and noting that  $(qr)^{s_n}/k_n \to \infty$, we have 
\begin{equation}\rho_n \to \rho.\label{eq:rhon}\end{equation}

For a set of vertices $A$ in the graph $\hat G_n$, let $y^{(0)}(A) = A,\; y^{(1)}(A)= y(A) = \{y_i(x): x\in A,\;1\leq i \leq r\}$ and $y^{(k+1)}(A) = y(y^{(k)}(A))$ for $k \geq 0$. Given a vertex $x$ and $R \in \N$, we define the ball $B(x, R)$ as the subgraph of $\hat G_n$ with vertices $\cup_{k=0}^R \; y^{(k)}(x)$ and all the edges of $\hat E_n$ that start and end at these vertices. Let $F(x,R)$ denote the event that $B(x,R)$ has no cycles and $F(x,y,R)$ the event that $B(x,R)$ and $B(y,R)$ have no cycles and are disjoint. Revealing edges one by one, it is easy to check that, by the choice of $s_n$,
\begin{equation}\label{eq:Fn}\lim_{n\to \infty}\P_n(F(x_1, s_n))=\lim_{n\to \infty}\P_n(F(x_1, x_2, s_n)) = 1.\end{equation}
If $F(x_1, s_n)$ occurs, then $(|\hat \xi^{x_1}_t|)_{0\leq t \leq s_n}$ has the same distribution as $(Z_t)_{0 \leq t \leq s_n}$ and, if $F(x_1, x_2, s_n)$ occurs, then $(|\hat \xi^{x_1}_t|)_{0\leq t \leq s_n},\;(|\hat \xi^{x_2}_t|)_{0\leq t \leq s_n} $ are distributed as two independent copies of this process.

For fixed $n$, let $X_i = I_{\left\{|\hat \xi^{x_i}_{s_n}| > k_n\right\}}$. Under $\P_n$, the random vector $(X_1,\ldots, X_n)$ is exchangeable and
$$|\E_n(X_1) - \rho_{n}| \leq \P_n(F(x_1, s_n)^c),\qquad \text{Cov}(X_1, X_2) \leq \P_n(F(x_1, x_2, s_n)^c).$$
Now, (\ref{eq:rhon}) and (\ref{eq:Fn}) imply that $\frac{1}{n}\sum_{i=1}^n X_i$ converges to $\rho$ in probability, as desired.
\end{proof}

\begin{proposition}\label{prop:afsn}There exists $c > 0$ such that
$$\lim_{n \to \infty}\P_n\left(\sup_{A \subset V_n: |A| \geq k_n} \;P_{\hat G_n}\left(\hat \xi^A_{e^{cn}} = \varnothing \right) \leq n^{-\sqrt{\log n}} \right) = 1.$$
\end{proposition}
Proving this result takes most of our effort. We postpone the proof and first show how the two propositions are used to establish the main theorem.\medskip\\

\bprthm{thm:main}
Let $c$ be the constant of Proposition \ref{prop:afsn}. Fix $\epsilon > 0$ and a sequence $(t_n)$ as in the statement of the theorem. If $t \leq e^{cn}$,  we have
\begin{eqnarray}
&&\nonumber\P_n\left(\frac{|\{x \in V_n: \hat \xi^x_t \neq \varnothing\}|}{n} \leq \rho - \epsilon \right) \leq \P_n\left(\frac{|\{x \in V_n: \hat \xi^x_{e^{cn}} \neq \varnothing\}|}{n} \leq \rho - \epsilon \right)\\
&&\leq \P_n\left(\frac{|\{x \in V_n: |\hat \xi^x_{s_n}| > k_n\}|}{n} \leq \rho - \epsilon\right) + \P_n\left(\exists x\in V_n: |\hat \xi^x_{s_n}|>k_n,\; \hat \xi^x_{e^{cn}} = \varnothing \right).\label{eq:auxThm}
\end{eqnarray}
The first term vanishes as $n \to \infty$ by Proposition \ref{prop:bfsn}. Let $H$ be the event inside $\P_n(\cdot)$ in the statement of Proposition \ref{prop:afsn}. The second term in (\ref{eq:auxThm}) is less than
$$\P_n(H^c) + n\cdot \frac{1}{n^{\sqrt{\log n}}} \stackrel{n \to \infty}{\xrightarrow{\hspace*{0.8cm}}} 0.$$
This shows that 
\begin{equation}
\label{eq:conct1}\inf_{t \leq e^{cn}}\P_n\left(\frac{|\{x \in V_n: \hat \xi^x_t \neq \varnothing\}|}{n} > \rho - \epsilon \right)\stackrel{n \to \infty} {\xrightarrow{\hspace*{0.8cm}}} 1.
\end{equation}

Now let us consider the reverse inequality. If $t \geq t_n$, we have
$$\P_n\left(\frac{|\{x \in V_n: \hat \xi^x_t \neq \varnothing\}|}{n} \geq \rho + \epsilon \right) \leq \P_n\left(\frac{|\{x \in V_n: \hat \xi^x_{t_n} \neq \varnothing\}|}{n} \geq \rho + \epsilon \right).$$
By a simplified version of the argument that established Proposition \ref{prop:bfsn}, it can be shown that the right-hand side vanishes as $n \to \infty$. Thus,
\begin{equation}
\label{eq:conct2}\inf_{t \leq e^{cn}}\P_n\left(\frac{|\{x \in V_n: \hat \xi^x_t \neq \varnothing\}|}{n} < \rho + \epsilon \right)\stackrel{n \to \infty} {\xrightarrow{\hspace*{0.8cm}}} 1.
\end{equation}

By (\ref{eq:conct1}), (\ref{eq:conct2}) and duality, we get
$$\inf_{t_n\leq t \leq e^{cn}}\P_n\left(\rho - \epsilon < \frac{|\xi^{V_n}_t|}{n} < \rho+\epsilon \right) \stackrel{n \to \infty}{\xrightarrow{\hspace*{0.8cm}}} 1$$
as required.
\eprthm

We now need to prove Proposition \ref{prop:afsn}; three preliminary results will be needed: Lemmas \ref{lem5domber}, \ref{lem:robgood} and \ref{lem:robeasy}. 

Once and for all, fix $\tilde q < q,\;\delta > 0$ and $g \in \mathbb{N}$ so that 
$$\tilde q r > 1, \;\delta < \min\left((\tilde q r -1), 1\right)\text{ and } (\tilde q r - 1- \delta)(\tilde q r)^{g-1} > 1+\delta.$$
We now give some definitions and notation.

Given $m \in \N$, let
\begin{eqnarray*} T^0_m &=& \{1,\ldots, m\},\\T^i_m &=& \{1, \ldots, m\} \times \{1, \ldots r\}^i, \quad 1 \leq i \leq g,\\
T_m &=& \cup_{i=0}^g\; T_m^i.\end{eqnarray*}
For $\sigma = (\sigma_0, \ldots, \sigma_i), \sigma ' = (\sigma '_0, \sigma '_1, \ldots, \sigma '_j) \in T_m$, we say $\sigma \prec \sigma '$ either if $i < j$ or if $i = j$ and $\sigma$ is less than $\sigma '$ in lexicographic order. With this order, we can take an increasing enumeration 
\begin{equation}T_m = \{\sigma^1, \ldots, \sigma^{(1+r+ \ldots +r^g)m}\}\label{eq:enum}\end{equation}
Then, $T_m^0 = \{\sigma^1, \ldots, \sigma^m\}$ and, for $i \geq 1$, $T^i_m = \{\sigma^{(1 + r + \ldots + r^{i-1})m + 1}, \ldots, \sigma^{(1 + r + \ldots + r^i)m}\}$.

Next, we endow $T_m$ with directed edges by setting
\begin{equation*}\sigma \to \sigma ' \text{ if and only if } \sigma = (\sigma_0, \ldots, \sigma_i),\; \sigma ' = (\sigma_0, \ldots, \sigma_i, \sigma '_{i+1}) \text{ for some } i.\end{equation*} $T_m$ is thus the disjoint union of $m$ rooted, directed trees, each with $g$ generations above the root. If we can go from $\sigma$ to $\sigma'$ by following a path of oriented edges of the tree, we say that $\sigma$ is an \textit{ancestor} of $\sigma'$.

The set $\{0, 1\}^{T_m}$ will be called the space of configurations. Given vertex $\sigma \in T_m$ and configuration $\psi \in \{0, 1\}^{T_m}$, $\psi(\sigma) \in \{0, 1\}$ will denote the value of $\psi$ at $\sigma$.

Now assume $\hat G_n = (V_n, \hat E_n)$ is given and $A \subset V_n$ with $|A| = m$. We can enumerate $A = \{x_{j_1}, \ldots x_{j_m}\}$ in the order of the indices of $V_n$. Given $\sigma = (\sigma_0, \ldots, \sigma_i) \in T_m$ with $i > 0$, let $z^{\sigma} = y_{\sigma_{i}}(y_{\sigma_{i-1}}(\cdots (y_{\sigma_1}(x_{j_{\sigma_0}}))\cdots))$. Finally, define
$$\mathcal{A}^\sigma = \{z^{\sigma '} \in T_m: \sigma ' \prec \sigma\}.$$
We now present an algorithm to construct a configuration $\psi = \psi(A) \in \{0, 1\}^{T_m}$ from $A$. The index $j$ in the algorithm follows the enumeration given in (\ref{eq:enum}).
$$\begin{aligned}&\textbf{for } j = 1 \textbf{ to } m \\
&| \textbf{ set } \psi(\sigma^{j}) = 0;\\
&\textbf{for } j = m+1 \textbf{ to } (1 + r + \ldots + r^g)m\\
&\left|\begin{array}{l}\textbf{if } \left[\psi(\sigma) = 1 \text{ for some } \sigma \text{ ancestor of } \sigma^j\right] \textbf{ or } \left[z^{\sigma^{j}} \notin \mathcal{A}^{\sigma^j}\right]\\
\quad \textbf{ then } \text{set } \psi(\sigma^{j}) = 0\\
\quad \textbf{ else } \text{set } \psi(\sigma^{j}) = 1
\end{array}\right.
\end{aligned}$$

In words, vertices are inspected in order; the roots are all set to $0$ and the other vertices are set to 0 either if one of their ancestors has already been marked with a 1 or if their image under the map $\sigma \mapsto z^\sigma$  has never been seen before; otherwise they are set to 1. Figure 1 presents an example of the effect of the algorithm.

\begin{figure}[htb]
\begin{center}
\setlength\fboxsep{0pt}
\setlength\fboxrule{0pt}
\fbox{\includegraphics[width = 1.0\textwidth]{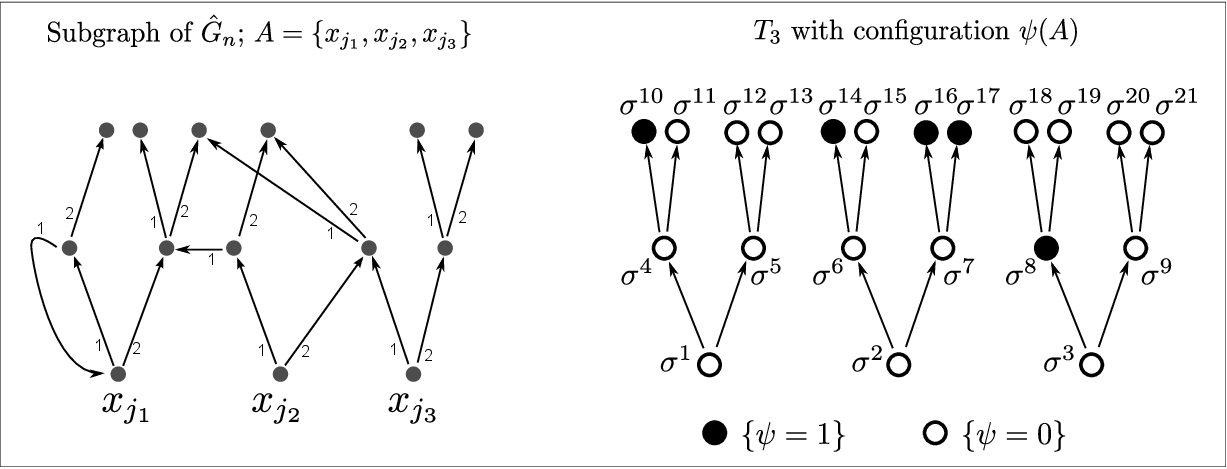}}
\end{center}
\caption{Example of the algorithm. Here $r = 2,\; g = 2$. The numbers in the arrows in the left diagram serve to distinguish $y_1(x)$ and $y_2(x)$ for each vertex $x$}
\end{figure}

\begin{lemma}
\label{lem5domber}
Given $A \subset V_n$ with $|A| = m$ and $\sigma^{i_1}, \ldots, \sigma^{i_k} \in T_m$,
$$\P_n\big(\;[\psi(A)](\sigma^{i_1}) = \ldots = [\psi(A)](\sigma^{i_k}) = 1\; \big) \leq \left(\frac{m + rm + \ldots + r^gm}{n - r}\right)^k.$$
\end{lemma}
\begin{proof} There is no loss of generality in assuming that $\sigma^{i_a} \prec \sigma^{i_b}$ when $a < b$. We then have
$$\begin{aligned}
&\P_n\big(\; [\psi(A)](\sigma^{i_k}) = 1 \;\big|\;[\psi(A)](\sigma^{i_1}) = \ldots = [\psi(A)](\sigma^{i_{k-1}}) = 1\; \big) \leq \frac{m + rm + \ldots + r^gm}{n - r}.
\end{aligned}$$
Indeed, let $\Theta^{i_k}$ denote the event that none of the ancestors of $\sigma^{i_k}$ in $T_m$ is marked with a $1$ in $\psi(A)$. First note that $\{[\psi(A)](\sigma^{i_k}) = 1\} \subset \Theta^{i_k}$, because the algorithm fills all positions above a 1 with 0's. Next, fix $a_{m+1}, a_{m+2}, \ldots, a_{i_{k}-1} \in V_n$ such that 
$$\begin{aligned}
&\{z^{\sigma^{m+1}} = a_{m+1}, \ldots, z^{\sigma^{i_{k}-1}} = a_{i_{k}-1}\} \subset \Theta^{i_k} \cap \{\; [\psi(A)](\sigma^{i_1}) = \ldots = [\psi(A)](\sigma^{i_{k}-1}) = 1\;\}
\end{aligned}$$
(we start at $m+1$ because $z^{\sigma^1}, \ldots, z^{\sigma^m}$ are always equal to the points of $A$). Then, conditioned on $\{z^{\sigma^{m+1}} = a_{m+1}, \ldots, z^{\sigma^{i_{k}-1}} = a_{i_{k}-1}\}$, there are at least $n-r$ possible positions for $z^{\sigma^{i_k}}$, and $[\psi(A)](z^{\sigma^{i_k}}) = 1$ precisely when $z^{\sigma^{i_k}} \in \mathcal{A}^{\sigma^{i_k}}$, a set of size less than $m + rm + \ldots + r^gm$.
\end{proof}

Given $A \subset V_n$ with $|A| = m$, let
$$d_i(A) = |\{\sigma \in T^i_m: [\psi(A)](\sigma) = 1\}|,\qquad d(A) = \sum_{i=1}^g d_i.$$
We say that $A$ is \textit{expansive} if $d(A) \leq (1+\delta)m$. The next lemma shows the motivation for this definition.

\begin{lemma}
\label{lem:robgood}
There exists $c_1 > 0$ such that, if $A \subset V_n$ is expansive, then
$$P_{\hat G_n}\left(|\hat \xi^A_{g}| < (1+\delta)|A| \right) \leq e^{-c_1|A|}.$$
\end{lemma}
\begin{proof}
Let $m = |A|$. If $i < g$ and $B \subset T^i_m$, we will write
$$J(B) = \{\sigma'\in T_m: \sigma \to \sigma' \text{ for some } \sigma \in B\}\subset T^{i+1}_m.$$
Consider the process $(\hat \xi^A_t)_{0 \leq t \leq g}$; define the sets
$$\begin{aligned}
&B_0 = \{\sigma \in T^0_m: z^\sigma \text{ gives birth at time } 0\};\\
&B_{i} = \{\sigma \in J(B_{i-1})\cap \{\psi(A)=0\}: z^\sigma \text{ gives birth at time }i\},\quad 1\leq i < g
\end{aligned}$$
The definition of $B_0$ implies that $\hat \xi^A_1 \supset \{z^\sigma:\sigma \in J(B_0)\}$. From the construction of $\psi(A)$ we see that $\sigma \mapsto z^\sigma$ is injective on $J(B_0)\cap \{\psi(A) = 0\}$, so we have $|\hat \xi^A_1| \geq |J(B_0) \cap \{\psi(A) = 0\}|$. Iterating this argument we get
\begin{equation}
\label{eq:auxClaim} |\hat \xi^A_i| \geq |J(B_{i-1})\cap\{\psi(A) = 0\}|,\qquad 1 \leq i \leq g.
\end{equation}

Define the events
$$\begin{aligned}
&F_0 = \{|B_0| < \tilde q m \},\\
&F_i =  \{|B_i| < \tilde q\cdot |J(B_{i-1}) \cap \{\psi(A) = 0\}|\},\quad 1 \leq i < g.
\end{aligned}$$
We now claim that
\begin{equation}\label{eq:claim}\left({\mathop \cup_{i=0}^{g-1}}\;F_i\right)^c \subset \left\{|\hat \xi^A_g| \geq (1+\delta)|A|\right\}.\end{equation}

Indeed, if none of the $F_i$ occurs, we have
$$\begin{aligned}
&|B_0| \geq \tilde q m;\\
&|J(B_0) \cap \{\psi(A) = 0\}| \geq r\cdot |B_0| -d_1 \geq \tilde q r m - d_1;\\
&|B_1| \geq \tilde q \cdot|J(B_0) \cap \{\psi(A) = 0\}| \geq \tilde q^2 r m -\tilde q d_1;\\
&|J(B_1) \cap \{\psi(A) = 0\}| \geq r\cdot |B_1| - d_2 \geq (\tilde q r)^2m - \tilde q r d_1 - d_2;\\
&\cdots\\
&|J(B_{i-1})\cap \{\psi(A) = 0\}| \geq (\tilde q r)^{i}m - (\tilde qr)^{i-1}d_1 - (\tilde q r)^{i-2}d_2 - \cdots - \tilde q rd_{i-1} - d_{i}
\end{aligned}$$
for $i \leq g$. In particular, using $\tilde q r > 1$ and the definition of expansiveness, for $0 < i \leq g$ we have
\begin{equation} 
\label{eq:expanB}
|J(B_{i-1})\cap \{\psi(A) = 0\}| \geq (\tilde q r)^{i}m - (\tilde q r)^{i-1} d \geq (\tilde q r)^{i-1}(\tilde q r - 1 - \delta)m.
\end{equation}
By the choice of $g$, this gives
$$|J(B_{g-1})\cap \{\psi(A) = 0\}| \geq (1+\delta)m.$$
Together with (\ref{eq:auxClaim}), this proves (\ref{eq:claim}).

The proof of the lemma will thus be complete if we show that, for some $c_1 > 0$,
\begin{equation}\label{eq:missing}
P_{\hat G_n}\left({\mathop\cup_{i=0}^{g-1}}\;F_i \right) \leq e^{-c_1m}.
\end{equation}
We start by writing
$$P_{\hat G_n}\left({\mathop\cup_{i=0}^{g-1}}\;F_i \right) \leq P_{\hat G_n}(F_0) + \sum_{i=1}^{g-1}P_{\hat G_n}\left(F_i\;\left|\; {\mathop \cap_{j=0}^{i-1}}F_j^c\right.\right).$$
In order to bound the terms of this sum, we will need the estimate
\begin{equation}
\nonumber \P(\mathsf{Bin}(k,p) \leq xkp) \leq \exp\{-\upgamma(x)kp\} \quad \text{for all } x \in (0,1),
\end{equation}
where $\upgamma(x) = x\log x - x + 1$. This follows from Markov's inequality; see Lemma 2.3.3 in \cite{Dur}. We then have
$$P_{\hat G_n}\left(F_0 \right) = \P(\mathsf{Bin}(m, q) < \tilde qm) \leq \exp\{-\upgamma(\tilde q/q)qm\}$$
Also, on the event $\cap_{j=0}^{i-1} F_j^c$, by (\ref{eq:expanB}) we have $|J(B_{i-1})\cap\{\psi(A) = 0\}| > (\tilde q r - 1 - \delta)(\tilde q r)^{i-1}m > (\tilde q r - 1 - \delta)m$, so
$$P_{\hat G_n}\left(F_i\; \left|\;{\mathop\cap_{j=0}^{i-1}F_j^c} \right. \right) \leq \exp\{-\upgamma(\tilde q / q)q(\tilde q r - 1 - \delta)m\}.$$
The proof of (\ref{eq:missing}) is now complete.
\end{proof}

\begin{lemma}
\label{lem:robeasy}
There exists $\upkappa > 0$ such that, putting $K_n = \upkappa \cdot n$,
$$\P_n\left(\exists A \subset V_n: k_n \leq |A| \leq K_n,\;\psi(A) \text{ is not expansive} \right) \stackrel{n \to \infty}{\xrightarrow{\hspace*{0.8cm}}} 0.$$
\end{lemma}
\begin{proof}
For fixed $m$ we have
$$\begin{aligned}&\P_n\left(\exists A \subset V_n:|A| = m,\;\psi(A) \text{ is not expansive} \right) \leq \sum_{A:|A|=m}\P_n(\psi(A) \text{ is not expansive})\\
&\leq \sum_{A:|A|=m} \;\sum_{d= \lceil (1+\delta)m\rceil}^{(1+r+\cdots+r^g)m} \; \sum_{D \subset T_m:|D|=d} \;\P_n\left([\psi(A)](\sigma) = 1 \; \forall \sigma \in D \right).
\end{aligned}$$
We now bound $|\{D \subset T_m:|D|=d\}|$ by $2^{|T_m|}$ and use Lemma \ref{lem5domber} to bound the probability; the above is less than
$$\begin{aligned}
&{n \choose m}\;(1 + r + \cdots + r^g)\;m\; 2^{(1 + r + \cdots + r^g)m}\; \left(\frac{(1 + r + \cdots + r^g)m}{n - r}\right)^{(1+\delta)m}\\
&\leq \left(\frac{ne}{m}\right)^m C^m \left(\frac{m}{n} \right)^{(1+\delta)m}\left(\frac{n}{n-r} \right)^{(1+\delta)m} \leq \left(C\left(\frac{m}{n}\right)^\delta\right)^m;
\end{aligned}$$
here $C$ is a constant that only depends on $r, g$ and $\delta$, and whose value has changed in the last inequality. Now choose $\upkappa$ such that $C\upkappa^\delta < 1/e$. The probability in the statement of the lemma is then less than
$$\sum_{i=k_n}^{K_n}e^{-i} \leq \upkappa n e^{-(\log n)^2}\stackrel{n \to \infty}{\xrightarrow{\hspace*{0.8cm}}} 0.$$
\end{proof}

\bprprop{prop:afsn}
Assume that $n$ is large enough that $\delta k_n > 1$ and that $\hat G_n$ satisfies
\begin{equation}
\label{eq:frombef}
\text{for every } A \subset V_n \text{ with } k_n \leq |A| \leq K_n,\;\psi(A) \text{ is expansive}.
\end{equation}
Let $c = \frac{c_1\upkappa}{2}$, where $c_1$ and $\upkappa$ are the constants of the two previous lemmas. We will prove that
\begin{equation}
\label{eq:wantL}
\text{for every } A \subset V_n \text{ with } |A| \geq k_n,\;P_{\hat G_n}\left(\hat \xi^A_{e^{cn}} = \varnothing \right) < n^{-\sqrt{\log n}}.
\end{equation}
Together with Lemma \ref{lem:robeasy}, this will imply the result we need.

We start noting that, if $|A| > k_n$, then 
\begin{equation}\label{eq:finalF}P_{\hat G_n}\left(|\hat \xi^A_g| < \min\left(|A|+1, K_n \right)\right) < e^{-c_1\min(|A|,K_n)}.\end{equation} 
Indeed, if $|A| < K_n$, this follows directly from Lemma \ref{lem:robeasy} and $(1+\delta)|A| > |A| + \delta k_n > |A| + 1$. If $|A| \geq K_n$, we can take a subset $A'\subset A$ with $|A'| = \lfloor K_n \rfloor$ and use the previous argument for $A'$ together with the fact that $\xi^{A'}_g \subset \xi^A_g$.

Using (\ref{eq:finalF}), we have
$$\begin{aligned}&P_{\hat G_n}\left(|\hat \xi^A_{j\cdot g}| \geq \min(|A| + j, K_n) \; \text{ for } 1 \leq j \leq e^{cn}\right) \geq 1 - \sum_{j=0}^{e^{cn}} e^{-c_1\min(|A|+j, K_n)} \\
&\qquad\qquad\qquad \qquad\qquad \qquad\geq1-\sum_{j=0}^{\lfloor K_n - k_n \rfloor} e^{-c_1(k_n+j)} - \sum_{j=\lfloor K_n - k_n \rfloor+1}^{e^{cn}} e^{-c_1 K_n} \\&\qquad \qquad\qquad \qquad\qquad \qquad\geq 1 - K_n \cdot e^{-c_1k_n} - e^{cn} \cdot e^{-c_1K_n}
\\&\qquad \qquad\qquad \qquad\qquad \qquad\geq 1 - \upkappa ne^{-c_1(\log n)^2} -e^{-\frac{c_1\upkappa n}{2}} > 1 -n^{-\sqrt{\log n}}
\end{aligned}$$
when $n$ is large enough, proving (\ref{eq:wantL}).
\eprprop


\end{document}